\documentclass[11pt]{article}
\usepackage{amssymb,amsmath,amsthm,amsfonts,fullpage}
\usepackage[dvips]{graphicx,color}
\usepackage{epsfig}
\usepackage[mathscr]{eucal}

\newtheorem{defi}{Definition}[section]
\newtheorem{theo}{Theorem}[section]
\newtheorem{prop}[theo]{Proposition}

\newtheorem{coro}[theo]{Corollary}

\begin{document}

\title{Topologically Fragmental Space and the Proof of Hadwiger's Conjecture}
\author{Cao Zexin
\thanks{
e-mail: caozx@mail.syiae.edu.cn}\\
Department of Mathematics and Physics
\\ Shenyang Institute of Aeronautical Engineering\\
110034 China }

\maketitle \setcounter{page}{1}
\begin{abstract}
A topological space is introduced in this paper. Just liking the
plane, it's continuous, however its $n+1$ regions couldn't be
mutually adjacent. Some important phenomenon about its
cross-section are discussed. The geometric generating element of
the coloring region-map is also an important concept. Every
$n$-coloring region map is in the cross-section set of an
$n$-color geometric generating element. The proof of four color
theorem and Hadwiger's conjecture is obtained by researching them
and their cross-sections. And we can see in the context, that
those conjectures are not of graph theory, but such topological
space.

{\bf keywords:}
topologically fragmental space, adjacent, generating element,
cross-section

\end{abstract}

\section{Introduce}

The Original Four Color Theorem is one of the simplest
mathematical problems to state and understand\cite{for0}. It says
that any maps are four-colorable. In other words, one can color
the map's regions in four colors at most , such that no two
adjacent regions are colored by the same color. The proof aided by
computer of the four color theorem was published by Appel and
Haken in 1977\cite{for1,for2,for3,for4}. A new proof following
this way is given by N. Robertson, D. P. Sanders, P. D. Seymour
and R. Thomas recently\cite{for5}. However, those proof aided by
computer are very different from here, and, I just think, the idea
in those proof might not be helpful to the proof of Hadwiger's
Conjecture for them depending on the planarity so deeply.

The coloring of geographical maps is essentially a topological
problem, in the sense that it only depends on the connectivities
between the countries, not on their specific shapes, sizes, or
positions. We can just as well represent every country by a single
point (vertex), and the adjacency between two bordering countries
can be represented by a line (edge) connecting those two points.
Then the dual graph could be taken as the equivalent of the planar
map. However, when the planar map is replaced by its dual graph,
its geometric property is neglected. For instance, the map could
fill the whole plane, the dual graph can't; and cutting the map,
its section, a 2-coloring line could be obtained, the product of
its dual graph is a few points. The most important is, a planar
map could be cut from a solid entity, its dual can't. And noticing
on the original four color conjecture\cite{for0}, it is about the
map, the plane, not the planar graph. That means, it might be the
plane that own the planar discreteness which is just inherited by
the planar graph. So the thing we wanted in those conjecture might
be geometric or topological. If it were, finding a proof in graph
theory might be a wrong way at the beginning.

The natural generalization of four color theorem was proposed by
Hadwiger in 1943\cite{had1}: (Hadwiger's conjecture) For all $n$,
every graph not contractible to $K^n$ is $(n-1)$-colorable. It was
represented in graph theory entirely. Hadwiger's conjecture looks
like more complicated than four color theorem and was also tried
by many people\cite{had2,had3,had4,had5}. This paper is just
another one. However, excepting this paper, there hasn't a visible
approach to this conjecture evenly, for the method used in proving
four color theorem no referring to it.

In the next, a region concept and its complete graph will be
introduced in graph theory firstly. Then the region is redefined
in the Euclidean space. Some other concepts such as map, complete
graph $K^n$ are also generalized, and their geometric meaning will
be shown. A topological space is defined on the regions and their
adjacent relationship. The important property of their
cross-section is also introduced. Then a proof of Four Color
Theorem and Hadwiger's conjecture is obtained by researching the
relation among them.

\section{The Region and Its Adjacency in Graph Theory} 
Firstly, some concepts are defined here.

\begin{defi}
A region $Y$ is a connected induced subgraph in graph $G$.
The minimal region is a vertex.
\end{defi}

For indicating its geometric meaning, a connected induced subgraph
is called as a region here. The region in a map is made up of
adjacent countries or one country. For its connectivity, every
region is contractible to a vertex obviously. The next is a
definition about the region's adjacent relationship.

\begin{defi}
There have two regions $Y_1, Y_2 \subset G$, $Y_1\cap Y_2 =
\emptyset$, $\forall a \in V(Y_1)$, $b \in V(Y_2)$, $\exists ab
\in E(G)$, then we say that the two regions $Y_1$, $Y_2$ are
adjacent, denoted as $Y_1|Y_2$.
\end{defi}

Two adjacent regions are contractible to two adjacent vertices,
respectively. By defining the adjacent relationship on the
subgraph, we can analyze graphs in a different way.

\begin{defi}
In graph $G$, $\exists Y_i, Y_i|Y_j$,$i\neq j, (i,j =
1,2,\cdots,n)$, we say the regions $Y_i(i=1,2,\cdots,n)$ make up a
complete $n$-region subgraph, denoted as $C^n$.
\end{defi}

Obviously, $C^n$ is contractible to a complete $n$-graph $K^n$.
For instance, It is easily to see that the Petersen
graph\cite{pet1} is a complete 5-region graph. Contracting its
every region to a vertex, a $K^5$ is obtained. In graph theory,
$C^n$ is equal to $K^n$ minor. However, the thing expressed by
$C^n$ is actually the macrostructure and the global character of a
graph. A graph could be split into five regions without thinking
one vertex to another and changing the graph. If those regions are
mutually adjacent, a $K^5$ minor should be in the graph. When
$K^n$ minor is replaced with $C^n$, the global property of graph
$G$ becomes clearer. And we can see in the next, its geometric
meaning also becomes clearer in this way. This next concept is
defined for describing the proof conveniently.

\begin{defi}
If graph $G$ is not contractible to $K^{n+1}$, but $K^n$, then we
say the complete number of $G$ is $n$. And the maximal complete
$n$-region subgraph in $G$ is $C^n$.
\end{defi}

Although it is important to the planarity, the bipartite graph
$K^{3,3}$ condition\cite{die2} is irrelative to the planar
coloring. Neglecting it, four color theorem and Hadwiger's
conjecture could be represented in the following form by using the
region.

\begin{theo}\label{t1}
Four Color Theorem: Every graph $G$ without a complete $5$-region
subgraph $C^5$ is $4-$colorable.
\end{theo}

\begin{theo}\label{t2}
Hadwiger's conjecture: Every graph $G$ without a $C^{n+1}$ is
$n-$colorable.
\end{theo}

\section{The Region and its Cross-section in Euclidean Space} 
Now we start inspecting the geometric aspect of those concepts.
This researching is fallen into an interesting mathematical field
which is called as region geometry\cite{geo1,geo2,geo3,geo4}. The
spatial relations such as adjacency, the duality are also an
important topic interested by computer geometer
specially\cite{adj1,adj2,adj3,adj4}. Here we just think that those
concepts are self-sufficient, and give some intuitive definition
about them.

As we know, a vertex in a planar graph is corresponding to a
country in the planar map. Here we can see, a country is a planar
region, and some connected countries make up a large region such
as a continent. Then we can give a definition of the planar
region, it is a non-empty connected open set with its boundaries.
When two planar regions are adjacent, they share a common boundary
segment, not just a point. Here, a vertex in graph theory is
corresponding to a planar region. Then what is dual to a vertex in
the nonplanar graph or $K^5$. The planar regions can fill the
plane completely, but it can't fill the whole solid space. So it
should be a solid object dual to a vertex of a nonplanar graph. We
call this object as a solid region. For instance, a piece of
plasticine is simply a solid region. When the region is solid, its
adjacent relationship should be different from the planar case. A
intuitive definition is following.
\begin{defi}
if two solid regions are coplanar, they are adjacent.
\end{defi}
A generalized definition of the region and its adjacent
relationship in $m$ dimensional Euclidean space is also following.
\begin{defi} \label{def32}
An $m$ dimensional region is an $m$-dimensional non-empty
connected open set with it boundaries in $\mathbb R^n$, and two
$m$ dimensional regions sharing an $m-1$ dimensional non-empty
region are adjacent.
\end{defi}
It need be emphasized that two $m$ dimensional regions only
sharing $m-2$ dimensional region are not adjacent, otherwise the
boolean relationship which graph theory is based on should be
broken as shown in Fig.1.
\begin{figure}[!ht]
\centering\epsfig{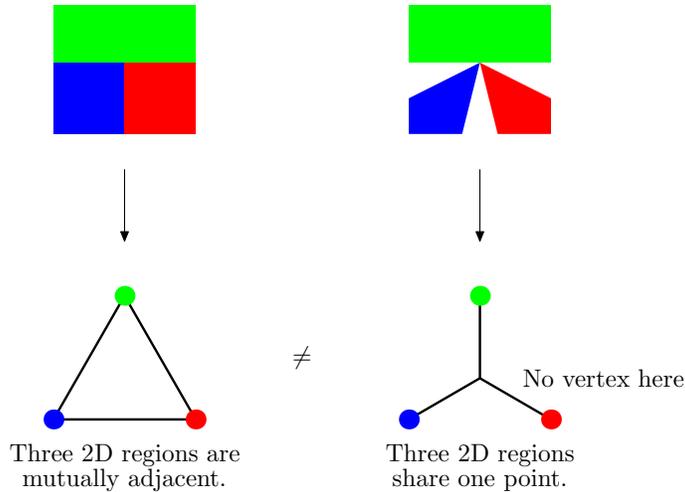}
\caption{Adjacency} \label{f1}
\end{figure}

Now we can consider any dimensional representations for graphs. It
is easy to see that any graph could be represented in Euclidean
space whose dimension is larger than 2 because any graph $G$ can
be embedded in it. We call the reality of graphs in $m$
dimensional Euclidean space as $m$ dimensional maps or $m$D-maps.
A $m$ dimensional map could be obtained by splitting $m$
dimensional Euclidean space into any mounts of $m$ dimensional
regions. The graph corresponding to a $m$D map is unchanged when
the map is under topological transformation, respectively.

\begin{figure}[!ht]
\centering\epsfig{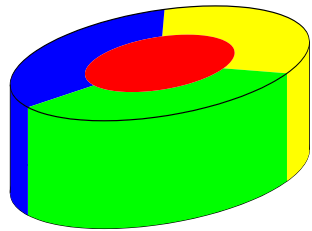} \caption{A solid
$K^4$ in $\mathbb R^3$} \label{f2}
\end{figure}

The $m$ dimensional reality of complete graph $K^n$ is introduced
here. For any $n$, complete graph $K^n$ can be embedded into
$m(m\geq 3)$ dimensional Euclidean space $\mathbb R^m$. Then we
just take every vertex as a center to make a corresponding little
$m$D region. The intersection of any two regions is null. Then
making every region extending along the edges which the vertex is
incident with until the two regions joined by the edge are
adjacent, an $m$ dimensional $K^n$ is obtained. A solid $K^n$ is
simply $n$ pieces of plasticine(regions) mutually adjacent. For
instance, a solid $K^4$ is shown in Fig.\ref{f2}.

It seems that no fresh thing appears by employing such a
representation. However, as mentioned above, the cross-sections of
a solid map are interesting. An important theorem is obtained soon
by researching the relationship between the solid $K^n$ and its
planar cross-sections.
\begin{theo}\label{t3}
Every $n$-colorable planar map is a cross-section of a solid
$K^n$.
\end{theo}
\begin{proof} Extending all the planar regions and their boundaries
along $Z$ axis a finite length $\delta (\delta>0)$, a solid map
whose adjacency is as same as the planar one is obtained. By
joining the same coloring regions with some thin plasticene poles
in the same color, those planar regions in same color are
integrated as one solid region now. There just have $n$ regions
mutually adjacent. So a solid $K^n$ is obtained.
\end{proof}

In the same way, a generalized corollary comes too.
\begin{coro}\label{t4} Every $m(m \geq 3)$ dimensional $n$-colorable map is a
cross-section of an $m+1$ dimensional $K^n$.
\end{coro}
Comparing to theorem \ref{t3} and \ref{t4}, the following
 are obvious.
\begin{coro}\label{tt3}
No planar map whose chromatic number is $(n+1)$ is in the
cross-section set of a solid $K^n$.
\end{coro}
\begin{coro}\label{tt4} No $m(m \geq 3)$ dimensional $n$-coloring map is
in the cross-section set of an $m+1$ dimensional $K^n$.
\end{coro}
It says in Theorem {\ref{t3}} and Corollary {\ref{t4}}, that all
the cross-sections of an $m+1$ dimensional $K^n$ make up a
complete set of $m$ dimensional $n$-coloring map. For instance, a
5-coloring planar map can't be cut from a solid $K^4$ absolutely
because all the cross-sections of a solid $K^4$ are 4-colorable.
The geometric meaning of the solid $K^n$ comes out, it is actually
\textbf{the visual minimal generating element} of all
$n$-colorable planar maps. And in general, an $m+1$ dimensional
$K^n$ is the visual minimal generating element of all the $m$
dimensional $n$-colorable region maps. So, as we see, the relation
between $K^n$ and $n$-coloring graph is easily represented in the
geometric way. A proof of four color theorem seems, if there had
no a solid $K^5$ in $\mathbb R^3$, it should be sure that the
plane is 4-colorable. But the solid $K^5$ does exist in $\mathbb
R^3$, so if a proof is wanted, we can try to construct a special
3D space which the solid $K^5$ can't be embedded in.

\begin{figure}[!ht]
\centering\epsfig{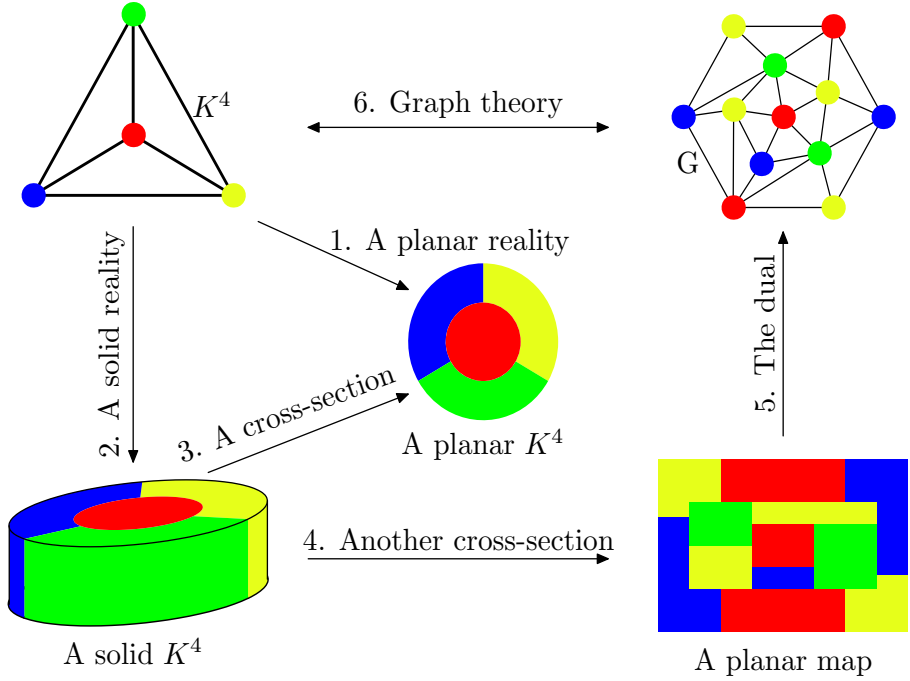} \caption{The
relation among $K^4$, a solid $K^4$, a planar $K^4$ and a planar
map} \label{f3}
\end{figure}

The relation among those concepts in the above is illustrated in
Fig.\ref{f3}. The relation (6) in Fig.\ref{f3} is researched in
graph theory, however, here the interesting is the relation (4).
Theorem {\ref{t3}} which can't be even depicted in graph theory,
is obtained easily by researching the relation (4). And some more
result are obtained and introduced later.

It need be mentioned that the meaning of the geometric $K^n$ is
more than it shown in the above. Different from a vertex, the
geometric regions can be split into more regions whose dual is a
connected induced subgraph in graph theory. So, the object
corresponding to a geometric region is actually a region in graph
theory, not a vertex. Thus the geometric $K^n$ is actually dual to
a complete $n$-region graph. To emphasize this feature, the
geometric $K^n$ is also called as the geometric $C^4$, and we
don't distinguish them later.


\section{Topologically Fragmental Space}  

Besides the geometric $C^n$, the plane also shows some important
property in this viewpoint. The plane is nothing but the complete
set of planar graph in graph theory. Now its topological feature
appears. Before discussing it, we can take a look at a simpler
case, the 2-colorable line $\mathbb R^1$.

According to the definition \ref{def32}, a region in the line is
just a segment. Two segments are adjacent when they share a point.
As everyone knows, no three segments are mutually adjacent in the
line, and no four segments are mutually adjacent in the cycle.
Such a feature shows that both the line and the cycle is discrete
somewhere, however, the line is well-known as one dimensional {\bf
continuous} Euclidean space $\mathbb R^1$. Comparing to the line,
the plane also has the similarly topological feature that no five
planar regions are mutually adjacent in it. So as we find, no $n$
regions mutually adjacent is the same feature of those space. A
sequacious supposal is that this feature is about a kind of
strange topological space. Although arbitrary solid $C^n$ could be
embedded into $\mathbb R^3$ that seems no such property, a special
solid space with the similar discrete attribution could be also
constructed in it. Thus, the topological space could be defined as
the following in general.

\begin{defi} Splitting $m$ dimensional Euclidean space $R^m$ into any
amount of $m$ dimensional regions(any adjacent regions can be
taken as one region), the set in which doesn't exist an $m$
dimensional complete regions-map $C^{n+1}$ is called as an $m$
dimensional topologically fragmental space, denoted as $C^nR^m$.
$C^n$ is an operator here, and $n$ is the complete number of this
space.
\end{defi}

Now it is clear that the line is $C^2\mathbb R^1$ for no three
segment mutually adjacent in it, and the plane $\mathbb R^2$ is
actually $C^4\mathbb R^2$. We can also construct a 3D fragmental
space $C^4\mathbb R^3$ in $\mathbb R^3$. $C^4\mathbb R^3$ has the
similar discrete structure to the plane $C^4\mathbb R^2$,
excluding that the solid object dual to the bipartite graph
$K_{3,3}$ could be embedded in. Now we know, the operator $C^n$
changes the space structure. For instance, the solid $K^4$ shown
in fig.\ref{f2} can't be embedded in $C^4\Bbb R^3$. In $C^4\Bbb
R^3$, the forth solid region is enveloped by other regions.

Two relations are interesting now, one is the relation between
$C^n\mathbb R^m$ and $C^n\mathbb R^{m+1}$, another is the relation
between $C^n\mathbb R^m$ and $m+1$ dimensional $C^j$. By
researching the first one, a theorem comes immediately.

\begin{theo}\label{t5}
The plane is the proper cross-section subset of $C^4\mathbb R^3$.
\end{theo}

\begin{proof} Extending all planar countries and their boundaries along
$Z$ axis to $(-\infty, \infty)$, a solid regions map is obtained.
The planar map is obviously its cross-section. And the graph dual
to this solid map is the same as the planar one. Any planar map is
in $C^4\mathbb R^2$, so no solid $C^5$ is in the solid map too. So
this solid map is in $C^4\mathbb R^3$.
\end{proof}
Any planar map are in the section set of $C^4\mathbb R^3$. It
means that $C^4\mathbb R^3$ is enough to generate the plane. As a
generalization of Theorem \ref{t5}, we have the following
corollary.
\begin{coro}\label{t6} $C^n\mathbb R^m$ is the cross-section subset of $C^n\mathbb
R^{m+1}$, $C^n\mathbb R^m \subset C^n\mathbb
R^{m+1}|_{x_{m+1}=0}$.
\end{coro}
The space $C^n\mathbb R^m$ has a very strange and notable
property, its cross-section set is not generally its subspace for
their different discrete structure.

\begin{figure}[!ht]
\centering\epsfig{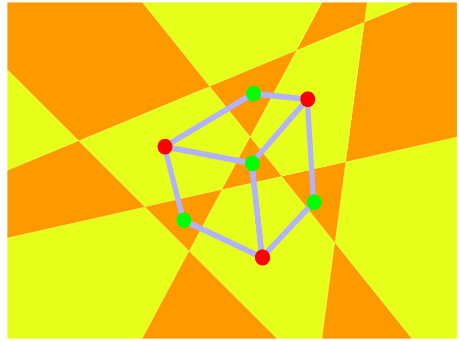} \caption{A
section map of a solid $K^2$ in $C^4\mathbb R^2$} \label{f4}
\end{figure}

\begin{figure}[!ht]
\centering\epsfig{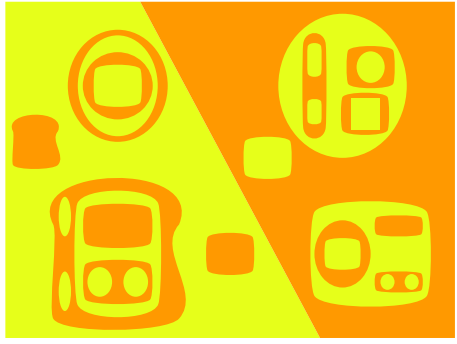} \caption{A
section map of a solid $K^2$ in $C^2\mathbb R^2$} \label{f5}
\end{figure}

Another interesting thing is the cross-section of a solid $C^n$ in
the plane, or the one of a $m+1$ dimensional $C^n$ in $C^n\mathbb
R^m$ in general. As we know, a planar $C^4$ could be embedded in
the plane, its cross-section of is just 2-colorable. And the
planar cross-section map of the solid $C^5$ seems to be
5-coloring, it is actually 4-colorable in a verifiable range
although rearranging its coloring is a very complicated operation.
However, such a phenomenon didn't happen on the cross-section of a
solid $C^i(i=2,3,4)$. Following the way in proving theorem
\ref{t3}, it is easily found that Fig.\ref{f4} is a cross-section
of a solid $C^2$. As Fig.\ref{f4} shows, a planar $C^4$ can be
even found in the section of a solid $C^2$. By analyzing those
case, we can obtain a conclusion that the complete number of a
cross-section is decided by the space $C^n\mathbb R^m$ that it is
located, to be more exact, decided by the operator $C^n$. For
another instance, whatever the solid $C^n(n>1)$ is, its
cross-section in the space $C^2\mathbb R^2$ is just the one as
Fig.\ref{f5} shows.

So, getting a cross-section of a solid $K^n$ have some choice now.
A planar cross-section in $C^2\mathbb R^2$, $C^3\mathbb R^2$ or
$C^4\mathbb R^2$? They are different. There have two steps in
getting a cross-section now. One operation is cutting the $m+1$
object in $m$ dimensional, denoted it as an operator $T$. For
instance, $T\mathbb R^{m+1} = \mathbb R^{m+1}|_{X_{m+1}=0} =
\mathbb R^m$. Another operation is $C^n$, no $C^{n+1}$ exists in
the cross-section map which it is operating on. The Corollary
\ref{t6} can be rewritten in the following.
\begin{coro}\label{tt6} $C^n T\mathbb R^{m+1} \subset T C^n\mathbb
R^{m+1}$.
\end{coro}
Corollary \ref{tt6} is an important theorem. When the two
operators operate on the solid space, it is shown in corollary
\ref{tt6} that what will happen on a solid $K^n$ while their
position are exchanged. For instance, considering the cutting
operator $T$ operating on a solid $K^5$, the result should be
5-coloring. After the operator $C^4$ operates on the result, a
planar map whose coloring isn't known is just obtained. We don't
know whether the operator $C^4$ changes the coloring of the
cross-section yet.

However, the coloring of the cross-section is clear when we
consider $TC^4(K^5)$. One among those region's adjacent
relationship is certainly broken when the solid $K^5$ is operated
by $C^4$. So before the cutting operator $T$ is operated on,
$C^4(K^5)$ has been already 4-colorable, then its cross-section is
4-colorable at most. The proof of four color theorem and
Hadwiger's Conjecture can be obtained now.

\section{ The Proof of Four Color Theorem and Hadwiger's Conjecture}
Before proving those conjecture, we represent them in geometric
form. Firstly, four color theorem is represented as the following.
\begin{prop}
Every planar map is a cross-section of a solid $K^4$.
\end{prop}
\begin{proof} We know that every $n$-coloring planar map is a
cross-section of a solid $K^n$ from theorem \ref{t3}.
If the chromatic number of a planar map is 5, it is impossible
that such a planar map is cut from a solid $K^4$. So if a
5-coloring planar map existed, a solid $K^5$ should be its
generating element at least.

However, from theorem \ref{t5}, every planar map is in the
cross-section set of $C^4\mathbb R^3$. According to Corollary
\ref{tt6}, the result that the solid $K^5$ is operated by $C^4 T$
is the subset of the one operated by $TC^4$. When this operator
$C^4$ operates on the solid $K^5$, one of the adjacent
relationship among $K^5$'s regions must be broken. So those
adjacent relationship in its cross-section corresponding to this
adjacent relationship is certainly broken too. It means, since no
solid $K^5$ exists in $C^4\mathbb R^3$, the cross-section map
corresponding to the solid $K^5$ specially, is forbidden in the
cross-section set of $C^4\mathbb R^3$. And according to theorem
\ref{t5} and \ref{t3}, so every planar map is a cross-section of a
solid $K^4$ and is 4-colorable.
\end{proof}
By using the same method, a corollary about Hadwiger's conjecture
is obtained easily.
\begin{coro}
Every solid map in $C^n\mathbb R^3$ is a cross-section of a
4-dimensional $C^n$.
\end{coro}

For understanding it better, the proof is demonstrated in
pictures. It is difficult to demonstrate the 3D case. I just
demonstrate the proposition
\textcolor[rgb]{0.00,0.00,1.00}{\textit {"the line is
2-colorable"}}, and the method proving it is the same as those
conjectures.
\begin{figure}[!ht]
\centering\epsfig{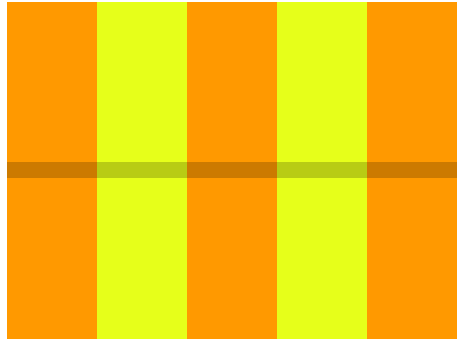} \caption{Every
linear `map' is a cross-section of $C^2\Bbb R^2$} \label{f6}
\end{figure}
The proof of theorem \ref{t5} is illustrated in Fig.\ref{f6}.
Taking the line as the axis $X$ in $\Bbb R^2$, then extending all
the segments along the axis $Y$ to $(-\infty, \infty)$, a planar
map in $C^2\mathbb R^2$ is obtained. As Fig.\ref{f6} shows, every
linear `map' is a section of $C^2\Bbb R^2$.
\begin{figure}[!ht]
\centering\epsfig{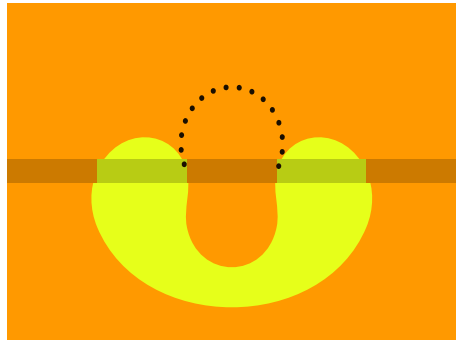} \caption{A
planar $K^2$ and the line} \label{f7}
\end{figure}
It is shown in Fig.\ref{f7} that every linear map is a
cross-section of a planar $K^2$. It is corresponding to theorem
\ref{t3}. And if taking the region surrounded by the yellow region
and dash line as the third region, the map is a planar $K^3$. When
it is operated by $C^2$, the third region is merged into the
orange region, or enveloped by the yellow region completely. So it
is very clear that after operated by $C^2$, its cross-section is
2-colorable.

Although theorem \ref{t3} and \ref{t5} are represented in this
approach easily and intuitively, it is difficult to represent them
in graph theory for their depending on the geometric feature yet.
However, as the proof shows, the intersection among topology,
geometry and graph theory such as topological graph
theory\cite{topg1}, might be a surprising mathematical field in
the future.

\end{document}